\documentclass{amsart}
\usepackage{amssymb,amsmath,
 theorem,euscript}

\newcounter{sec}

\newcounter{punct}[sec]

\begin{document}

\def\SL{\mathrm {SL}}
\def\SU{\mathrm {SU}}
\def\GL{\mathrm  {GL}}
\def\U{\mathrm  U}
\def\OO{\mathrm  O}
\def\Sp{\mathrm  {Sp}}
\def\SO{\mathrm  {SO}}
\def\SOS{\mathrm {SO}^*}

\def\PGL{\mathrm  {PGL}}
\def\PU{\mathrm {PU}}

\def\Gr{\mathrm{Gr}}

\def\Fl{\mathrm{Fl}}

\def\SpU{\mathrm {SpU}}
\def\OU{\mathrm {OU}}
\def\OGL{\mathrm{OGL}}

\def\GLO{\mathrm{GLO}}

\def\OSp{\mathrm {OSp}}

\def\Mat{\mathrm{Mat}}

\def\LLambda{\mathbf\Lambda}

\def\Pfaff{\mathrm {Pfaff}}

\def\B{\mathbf B}

\def\phi{\varphi}
\def\epsilon{\varepsilon}
\def\kappa{\varkappa}

\def\le{\leqslant}
\def\ge{\geqslant}

\renewcommand{\Re}{\mathop{\rm Re}\nolimits}

\renewcommand{\Im}{\mathop{\rm Im}\nolimits}

\def\spin{\mathrm{spin}}

\def\pia{\pi_\downarrow}

\def\graph{\mathrm{graph}}

\def\Diff{\mathrm{Diff}}

\newcommand{\im}{\mathop{\rm im}\nolimits}
\newcommand{\indef}{\mathop{\rm indef}\nolimits}
\newcommand{\dom}{\mathop{\rm dom}\nolimits}
\newcommand{\codim}{\mathop{\rm codim}\nolimits}

\def\cA{\mathcal A}
\def\cB{\mathcal B}
\def\cC{\mathcal C}
\def\cD{\mathcal D}
\def\cE{\mathcal E}
\def\cF{\mathcal F}
\def\cG{\mathcal G}
\def\cH{\mathcal H}
\def\cJ{\mathcal J}
\def\cI{\mathcal I}
\def\cK{\mathcal K}
\def\cL{\mathcal L}
\def\cM{\mathcal M}
\def\cN{\mathcal N}
\def\cO{\mathcal O}
\def\cP{\mathcal P}
\def\cQ{\mathcal Q}
\def\cR{\mathcal R}
\def\cS{\mathcal S}
\def\cT{\mathcal T}
\def\cU{\mathcal U}
\def\cV{\mathcal V}
\def\cW{\mathcal W}
\def\cX{\mathcal X}
\def\cY{\mathcal Y}
\def\cZ{\mathcal Z}

\def\frA{\mathfrak A}
\def\frB{\mathfrak B}
\def\frC{\mathfrak C}
\def\frD{\mathfrak D}
\def\frE{\mathfrak E}
\def\frF{\mathfrak F}
\def\frG{\mathfrak G}
\def\frH{\mathfrak H}
\def\frJ{\mathfrak J}
\def\frK{\mathfrak K}
\def\frL{\mathfrak L}
\def\frM{\mathfrak M}
\def\frN{\mathfrak N}
\def\frO{\mathfrak O}
\def\frP{\mathfrak P}
\def\frQ{\mathfrak Q}
\def\frR{\mathfrak R}
\def\frS{\mathfrak S}
\def\frT{\mathfrak T}
\def\frU{\mathfrak U}
\def\frV{\mathfrak V}
\def\frW{\mathfrak W}
\def\frX{\mathfrak X}
\def\frY{\mathfrak Y}
\def\frZ{\mathfrak Z}

\def\fra{\mathfrak a}
\def\frb{\mathfrak b}
\def\frc{\mathfrak c}
\def\frd{\mathfrak d}
\def\fre{\mathfrak e}
\def\frf{\mathfrak f}
\def\frg{\mathfrak g}
\def\frh{\mathfrak h}
\def\fri{\mathfrak i}
\def\frj{\mathfrak j}
\def\frk{\mathfrak k}
\def\frl{\mathfrak l}
\def\frm{\mathfrak m}
\def\frn{\mathfrak n}
\def\fro{\mathfrak o}
\def\frp{\mathfrak p}
\def\frq{\mathfrak q}
\def\frr{\mathfrak r}
\def\frs{\mathfrak s}
\def\frt{\mathfrak t}
\def\fru{\mathfrak u}
\def\frv{\mathfrak v}
\def\frw{\mathfrak w}
\def\frx{\mathfrak x}
\def\fry{\mathfrak y}
\def\frz{\mathfrak z}

\def\fros{\mathfrak{s}}

\def\bfa{\mathbf a}
\def\bfb{\mathbf b}
\def\bfc{\mathbf c}
\def\bfd{\mathbf d}
\def\bfe{\mathbf e}
\def\bff{\mathbf f}
\def\bfg{\mathbf g}
\def\bfh{\mathbf h}
\def\bfi{\mathbf i}
\def\bfj{\mathbf j}
\def\bfk{\mathbf k}
\def\bfl{\mathbf l}
\def\bfm{\mathbf m}
\def\bfn{\mathbf n}
\def\bfo{\mathbf o}
\def\bfp{\mathbf q}
\def\bfr{\mathbf r}
\def\bfs{\mathbf s}
\def\bft{\mathbf t}
\def\bfu{\mathbf u}
\def\bfv{\mathbf v}
\def\bfw{\mathbf w}
\def\bfx{\mathbf x}
\def\bfy{\mathbf y}
\def\bfz{\mathbf z}

\def\bfA{\mathbf A}
\def\bfB{\mathbf B}
\def\bfC{\mathbf C}
\def\bfD{\mathbf D}
\def\bfE{\mathbf E}
\def\bfF{\mathbf F}
\def\bfG{\mathbf G}
\def\bfH{\mathbf H}
\def\bfI{\mathbf I}
\def\bfJ{\mathbf J}
\def\bfK{\mathbf K}
\def\bfL{\mathbf L}
\def\bfM{\mathbf M}
\def\bfN{\mathbf N}
\def\bfO{\mathbf O}
\def\bfP{\mathbf P}
\def\bfQ{\mathbf Q}
\def\bfR{\mathbf R}
\def\bfS{\mathbf S}
\def\bfT{\mathbf T}
\def\bfU{\mathbf U}
\def\bfV{\mathbf V}
\def\bfW{\mathbf W}
\def\bfX{\mathbf X}
\def\bfY{\mathbf Y}
\def\bfZ{\mathbf Z}

\def\la{\langle}
\def\ra{\rangle}

\def\R {{\mathbb R }}
 \def\C {{\mathbb C }}
  \def\Z{{\mathbb Z}}
  \def\H{{\mathbb H}}
\def\K{{\mathbb K}}
\def\N{{\mathbb N}}
\def\Q{{\mathbb Q}}
\def\A{{\mathbb A}}

\def\T{\mathbb T}

\def\bbA{\mathbb A}
\def\bbB{\mathbb B}
\def\bbD{\mathbb D}
\def\bbE{\mathbb E}
\def\bbF{\mathbb F}
\def\bbG{\mathbb G}
\def\bbI{\mathbb I}
\def\bbJ{\mathbb J}
\def\bbL{\mathbb L}
\def\bbM{\mathbb M}
\def\bbN{\mathbb N}
\def\bbO{\mathbb O}
\def\bbP{\mathbb P}
\def\bbQ{\mathbb Q}
\def\bbS{\mathbb S}
\def\bbT{\mathbb T}
\def\bbU{\mathbb U}
\def\bbV{\mathbb V}
\def\bbW{\mathbb W}
\def\bbX{\mathbb X}
\def\bbY{\mathbb Y}

 \def\ov{\overline}
\def\wt{\widetilde}
\def\wh{\widehat}

\def\P{\mathbb P}

\def\arr{\rightrightarrows}

\def\SS{\smallskip}

\def\ev{{\mathrm{even}}}
\def\od{{\mathrm{odd}}}

\def\q{\quad}

\def\F{\mathbf F}

\def\b{\mathbf b}

\def\RA{\Longrightarrow}

\newfont{\weird}{cmff10}
\def\yu{\weird}
\def\yun{\weird YuN:\,\,}

\def\obmanka{$\mathop{}\quad$}

\title
{Mathematical researches of D.~P.~Zhelobenko}

\author
  {Yu.~A.~Neretin, S.~M.~Khoroshkin}

\begin{flushright}
\small Russian Mathematical Surveys, 2009, 1
\end{flushright}

\address{Yu.A.N.: University of Vienna,
Institute for Theoretical
and Experimental
Physics, MechMath Dept. of Moscow State University;
\newline neretin(at)mccme.ru,
\newline URL:
wwwth.itep.ru/$\sim$neretin, www.mat.univie.ac.at/$\sim$neretin}

\address{Khoroshkin S.M.:
Institute for Theoretical
and Experimental
Physics,
\newline
 khor(at)itep.ru}

\begin{abstract} This is a brief overview of  researches
of Dmitry Petrovich Zhe\-lo\-ben\-ko
(1934--2006).
He is the best known for his book "Compact Lie groups and their representations"
and for the classification of all irreducible representations of complex semisimple
Lie groups. We tell also on other his works, especially
on the spectral analysis of representations.
\end{abstract}



\maketitle

\bigskip

{\bf 1. Brief personal data.} The  mathematical researches of Dmitri\^{\i} Petrovich
Zhelobenko (Ul'yanovsk, 1934 - Moscow, 2006) are mainly devoted to the representation
theory of semisimple Lie groups and Lie algebras and to the noncommutative
harmonic analysis.

He  graduated from the Physical Department of
the Moscow State University (1958) and
took his postgraduate program  at  the Steklov Mathematical Institute under the supervision of
S.~V.Fomin and of M.~A.Naimark (1961). He has defended his PhD thesis
''{\it Harmonic Analysis
on the Lorentz group and some questions
 of the theory of linear representations}'' in 1962.
The doctoral thesis ''{\it Harmonic analysis of functions on semisimple Lie groups and its
applications to the representation theory}'' was defended in
 the Steklov Institute in 1972
 (official opponents: I.~M.~Gel'fand, I.~R.~Shafarevich, and A.~A.~Kirillov).
 In 1974 he was  an invited
 speaker at the Mathematical Congress in Vancouver%
 \footnote{In fact  he
was not let abroad at that time. By the way, Naimark's
 report at Nice Congress at 1970 was the joint
 with Zhelobenko by the contents.}

Since 1961 and for the rest oh his life
 he has worked in P.~Lumumba People's Friendship University%
\footnote{PFU, later renamed into People's Friendship University
of Russia (PFUR).}.

It appears that this job gave him an opportunity to use  his teaching abilities,
but there were not many not graduate students%
\footnote{besides he happened face certain administrative impediments when trying
to 'pick up' a gifted student. Here
 is  the list of postgraduates, who had defended PhD under
the supervision of Zhelobenko:
M.~S.~Al-Nator,
A.~S.~Kazarov, A.~G.~Knyazev, A.~V.~Lutsyuk.
We  took borrowed it
from his Curriculum Vitae. Partially, Vijay Jha
(he later worked in number theory) is also a student of Zhelobenko}.

In 1963--1978 M.~A.~Naimark, D.~P.~Zhelobenko, and A.~I.~Stern held a seminar
on representation theory in the Steklov Mathematical Institute. The lecturers were,
 among others:
 F.~A.~Berezin, N.~Ya.~Vilenkin, R.~S.~Ismagilov, A.~A.~Kirillov,
E.~V.~Kissin, G.~L.~Litvinov, V.~I.~Man'ko, M.~B.~Menskii, V.~F.~Molchanov,
M.~G.~Krein, G.~I.~Olshansky, I.~Sigal, V.~N.~Tolstoy, C.~Foia\c{s},
Dj.~Khadzhiev, P.~Halmos, A.~Ya.~Khelemsky, E.~Hewitt.
 Dmitri\^{\i} Petrovich gave lectures for physicists many times; they
took place in Dubna at 1965, in the College for Mathematical Physics of
the  Independent University of Moscow in early 90-ths,
 in the Sofia University in 1973--1974, and
on some physical conferences.

Since 1959 and
up to the last days of his life Dmitri\^{\i} Petrovich
was seriously ill (as a result of unlucky mountain hike)%
\footnote{Zhelobenko rescued a tourist group  in Altai mountains
in emergency conditions.
Here follows
a fragment of a letter of V.~I.~Chilin.
''I learned from the conversations with Dmitri\^{\i} Petrovich, that he got ill after
a hard mountain hike. After an agreement with the
leader of the group he went forward, searching the path. Then a heavy snowfall began
and he could not return back to the camp
 in the approaching dusk. He was 
was bound to circle  around
throughout the whole night
in order not to freeze. As a result he got cold developed later into a serious illness''.}.

Dmitri\^{\i} Petrovich is  the author of 9 books%
\footnote{Two of them are now in press, \cite{Zhe-P1}, \cite{Zhe-P2}.}
and of more than 70 papers (according to MathSciNet and Zentralblatt).
The goal of our short survey is to discuss the most significant, in our
regard,
 works of Zhelobenko. Clearly, this survey is not complete. Besides,
having set such a goal we then cannot deviate from the subject of
a purely professional text.
MCCME%
\footnote{The Moscow Center for Continuous Mathematical Education.}
 publishers and Yu.~N.~Torkhov are preparing now  edition of a new book by
Zhelobenko ''{\it Gauss algebras}''. It is assumed
to include  reminiscences of
his colleges and friends, G.~L.~Litvinov and
V.~F.Molchanov, and a paper of
his students  ,M.~S.~Al-Nator, I.~V.~Goldes,
A.~S.~Kazarov,  V.~R.~Nigmatullin.
The bibliography of his papers should be also present there.


\SS
{\bf 2. The book "Compact Lie groups and their representations" (1970).}
This is the
 classical book, which after 37 years remains one of the best books on
representation theory in both genres: as a textbook (or
an introductory book), and as
a monograph.
This text, peculiar in its structure as well as in its style%
\footnote{It seems that the evident popularity of his "wrongly" written
("for physicists") book surprised Dmitri\^{\i} Petrovich  up to the end of
 his life.}, finally became
  a specimen of "mathematics with a human face". It was designed for
 a reader with a background of 2-3 grades of a mathematical
 (or physical) department
 of that time;  it can serve as a textbook
 for the primary acquaintance or studying,
  as well as a book,
 interesting to specialists, or a handbook in a certain field.

We are not going to discuss the content (the book is worth reading itself), but note that
its "complicated chapters" (problems of spectral analysis) are in many respects based on
 the original results of D.~P.~Zhelobenko, in particular
 \cite{Zhe-spectral}, \cite{Zhe-lin}%
\footnote{It is worth to note that the paper by R.~Godement \cite{God}
influenced much on these works.}.

The original invention of Zhelobenko was the consistent exposition of a big number
of classical problems of the representation theory in the framework of so-called
"method of $Z$-invariants".
The method is based on a realization of
 representations of semisimple groups in
 functions on  the maximal unipotent
(strictly triangular) subgroup;
these functions must satisfy
 a certain system of partial
differential equations
\cite{Zhe-spectral}, \cite{Zhe-lin}. Virtually this method anticipates the theory of
$D$-modules on flag manifolds%
\footnote{Strictly triangular subgroup is a "big cell"
 of the flag variety}, see A.~Beilinson, I.~Bernstein \cite{BB},
J.-L.~Brylinsky, M.~Kashiwara \cite{BK}.

\SS


{\bf 3. The book "Harmonic Analysis on complex semisimple
Lie  groups" (1974).}  Essentially,
 this  is an original work written upon
 Zhelobenko's papers%
\footnote{Some of these works have overlaps with parallel
researches of M.~Duflo of that time, see references
in \cite{Duf}, \cite{Zhe-complex}.}
  \cite{Zh1}--\cite{Zh6}, \cite{Zh7}--\cite{Zh8},
   \cite{ZhN1}--\cite{ZhN2}  (1963--1973),
partially joint with  M.~A.~Naimark.
The main result is the Zhelobenko-Naimark theorem
on classification of all irreducible representations%
\footnote{in Frechet  spaces up to infinitesimal equivalence
of modules of $K$-finite vectors.}
of complex semisimple Lie groups. Let us recall this statement for
 the group
 $G=\GL(n,\C)$.%
\footnote{Below
$G$ is arbitrary semisimple group, $K$ is its maximal compact subgroup,
 $\frg\supset\frk$ are their Lie algebras. Now $G=\GL(n,\C)$,
 $K=\U(n)$.}

Denote by  $\U(n)$ the subgroup of all unitary matrices.
Denote by
$T\subset G$ the subgroup of upper-triangular matrices. For $A\in T$ denote by
$a_{jj}$ the diagonal elements of $A$. Fix a collection of complex numbers
 $p_1$,
\dots, $p_n$, $q_1$, \dots, $q_n$
 Consider the character  $T\to \C^*$, given by the relation
$$
\chi_{p,q}(A)=\prod_{j=1}^n a_{jj}^{p_j} \ov a_{jj}^{q_j}
.
$$
All  $(p_j-q_j)$ must be integers; otherwise  complex powers
make no sense.

Consider the representation $\rho_{p,q}$
 of the group $G$ induced from the representation
$\chi_{p,q}$ of the group $T$.%
\footnote{For the reader accustomed to another language. We consider
representations of $\GL(n,\C)$ in smooth sections of
complex linear bundles over a flag manifold.} What we get  is called
 the the {\it principal series representation}.
For generic  $p$, $q$%
\footnote{The exceptionally condition
 is $p_i-p_j\in\Z$, $q_i-q_j\in\Z$ for at least one pair
 $i$, $j$} principal series representations are irreducible. For exceptional
values of
 parameter the representation  $\rho_{p,q}$ has a finite
Jordan--H\"older
 composition series.

Next, for generic $p$,  $q$ a simultaneous
permutation of the tuples
$(p_1,\dots, p_n)$ and $(q_1,\dots,q_n)$ get an equivalent representation.%
\footnote{For the reader slightly familiar with the Lorentz group $\SL(2,\C)$ and
with the representation theory.  This
 is almost evident. First, for
 $\GL(2,\C)$ this is a well-known and easily checked statement.
Next, denote by
$T_{i}$ the generalized upper-triangular subgroup, where the element
 $a_{(i+1)i}$, which stands below the diagonal, is allowed to be nonzero.
  The Lie group  $T_{i}$ has $\SL(2,\C)$  as a semisimple factor.
Now we can construct our representation taking
a successive  induction,
first from
   $T$ to $T_{i}$, and   second
 to  $G$. But now for  $\SL(n,\C)$ we can permute
   the parameters  $(p_i,q_i)$ and $(p_{i+1},q_{i+1})$.} For exceptional values
   of $p$,
$q$ such representations are not equivalent,
but their composition series have the same
factors (we use the term "subfactor" below).

Now we formulate the Zhelobenko-Naimark theorem:
{\it  All irreducible representations of
$\GL(n,\C)$ are in one-to-one correspondence with
 orbits of the symmetric group  $S_n$ on the set of collections
$$(p_1,\dots,p_n\,; \,q_1,\dots,q_n)$$
Namely, for each orbit  one of the subfactors
of the corresponded
$\rho_{p,q}$  is taken.}

\SS

 Let us describe the subfactor more precisely.
Without loss of generality we can set
$$
p_1-q_1\ge p_2-q_2\ge \dots\ge p_n-q_n
.
$$
 Now one should take the subfactor (it is actually a subrepresentation),
 whose
  restriction to the subgroup
 $K=\U(n) $  contains a representation with the highest weight
 $(p_1-q_1,\dots,
p_n-q_n)$.

\SS

The theorem (1966) by itself
 was a part of thirty-year-long story around Harish-Chandra subfactor
theorem. Principal series and the parabolic induction were introduced in the
book  by M.~I.~Gel'fand and M.~A.~Naimark \cite{GN}. In 1954 Harish-Chandra
\cite{Harish} proved that any irreducible representation of a semisimple Lie
group can be realized as a subfactor of the principal series%
\footnote{see also \cite{Dix}, 9.4}.
 F.~A.~Berezin (1956), \cite{Ber0},
   \cite{Ber},
    has immediately made an attempt
 to turn this result into an explicit
classification in the case of complex groups.%
\footnote{Literally the paper yet fell short of the classification;
to be more precise the final theorem was given without  proof.
Serious gaps were found in this paper, and a discussion between Harish-Chandra
and Berezin on the correctness
of proofs followed.
  In the next paper ("{\it A Letter to Editors}") Berezin removed these gaps.
 Berezin's approach (bringing up to the classification)
  is expounded in  Appendix A of the discussed book of
 Zhelobenko
  (however we are bound to note that Zhelobenko
applied additional  arguments).
\newline
\obmanka
By the way, this work of Berezin contains a calculation of
radial parts of Laplace operators on semisimple groups, which was no less important
than the incipient "struggle" for the classification of representations.}
 Ten years later (1966) the Zhelobenko--Naimark Theorem has appeared,
  then in 1973 R.~Langlands \cite{Lan} provided a classification of all
  irreducible representations of semisimple groups, and Casselman
 \cite{Cas} in 1975 proved the subrepresentation theorem.%
\footnote{Any irreducible representation is a subrepresentation
in a principal series representation.}%
$^,$%
\footnote{The proofs were published in the book  by
A.~Borel and--N.~Wallach in 1980
\cite{BW} and in the paper by
B.~Casselman and D.~Mili{c}i{c} in 1982 \cite{CM}; see also the book  by
 A.~Knapp
\cite{Kna}. The real case in essence is more complicated compared to complex,
but the results are less transparent. The classification of all representations
is based on the classification of discrete series. Representations of discrete series
 themselves are rather complicated objects; a few descriptions of them are known,
 but  they all are still far from  transparent. After the work of
 M.~Flensted-Jensen\cite{F-J} had appeared it seemed that the discrete series
 were just about  to be quite understood.
 Thirty years afterwards the work of Flensted-Jensen remains unclear.}

\SS

Zhelobenko \cite{Zh4} also introduces the following structure.%
\footnote{This structure is  also
 considered in the book  by Dixmier
\cite{Dix} 1974,
Chapter 9, with a reference to the work of J.~Lepowsky and G.~W.~MacCollum,
this structure is quite natural and probably
was introduced before.}
Denote by
  $U(\frg)$ the enveloping algebra of the Lie algebra $\frg$,
 and by $U(\frk)$ the enveloping algebra of the subalgebra
$\frk$.
Denote by
 $\pi_\lambda$ the irreducible representation of the group $K$
 with the
highest weight $\lambda$.
By $\ker\pi_\lambda\subset U(\frk)$ we
denote the two sided ideal, which consists of  elements
 $v$, such that
$\pi_\lambda(v)=0$.  Next, by $U_{\lambda\mu}$ we denote the set of all
$u\in U(\frg)$ such that
$$(\ker\pi_\lambda) u\subset U(\frg)\, (\ker\pi_\mu)
.
$$
It is easy to check that
 $U_{\lambda\mu}U_{\mu\nu}\subset U_{\lambda\nu}$. Therefore we get a category
whose objects are dominant weights $\lambda$, and morphisms are elements
 $u\in U_{\lambda\mu}$.

For a representation of the group
  $G$ in the space
$H$ consider its restriction to the subgroup $K$; it decomposes into a direct
sum
$H=\oplus H_\lambda$, where in $H_\lambda$ the group $K$ acts by means of
representation $\pi_\lambda$ with certain multiplicity. It is easy to check that
 $u\in U_{\lambda\mu}$   sends $H_\mu$ to
$H_\lambda$, i.e, we get a representation of our category.

This subject, which is discussed in the book in detail,
can be naturally regarded as an
analog of the classical operational calculus.
The enveloping algebra is simultaneously the algebra of left-invariant
differential operators on the Lie group, and the action
of $U(\frg)$ in representations
is an analog of the Fourier transform of polynomial differential operators.

Integral formulas for intertwining operators for the principal nonunitary series
together with a detailed description of the adjoint action of the Lie algebra
 $\frk$ in $U(\frg)$ enable to describe the algebra
$U_{\lambda\lambda}/\mathrm{ker}\,\pi_\lambda U(\frg)$
 and classify all irreducible representations with the lowest weight $\lambda$.
 This result is then used in the classification
 of all irreducible representations of a complex semisimple Lie group.

\SS

From the  "sports" point of view,
the Zhelobenko-Naimark Theorem is left behind long ago.
On the other hand, the book still remains an instructive text
on complex semisimple groups (which are generally  understood better than real)
and together with accompanying papers could become a "basis" for  progressive
movement in the future. Unfortunately, the book is written in the style used in the
modern mathematics; it happens to be legally  precise text
for experts.%
\footnote{Afterwards Dmitri\^{\i} Petrovich regretted to have  submitted
to the influence
of "bourbakism", just coming into a fashion that time.}
As a result, the "operational calculus" is known and understood less than it deserves.
An important text on this subject is M.~Duflo's notes
\cite{Duf} 1974. \footnote{Under the influence of these notes, Zhelobenko's technique
 was used in the work of P.~Torasso \cite{Tor},
it was also used in one of the variants
 of the proofs of Baum-Connes conjecture by Vincent Lafforgue
 \cite{Laf}.}

Now, the unfinished Zhelobenko's book
"{\it Operational calculus on complex semisimple Lie groups}" is in press
(we did not see it yet).
Judging by the  title, it continues the subjects discussed above.
Hopefully it can help entering of these  topics
 to {
a common knowledge}.

\SS


{\bf 4. The book "Representations of reductive Lie algebras". 1994. }
This work and the related series of papers probably represent now the most
 interesting part of the original scientific heritage of
 D.~P.~Zhelobenko. The main part of the book is devoted to the exposition
 of the so-called reduction algebras, or Mickelsson algebras
 \cite{Mick}.

\SS

Let $\frg\supset\frk$  be complex Lie algebras, and
$\frk$  reductive.%
 \footnote{It is worth  to note the unexpected generality: subalgebra
 $\frk$ does not need not to be symmetric.} We are interested in the problem of
 restriction of finite-dimensional representations of
 $\frg$ to $\frk$.%
\footnote{It is reasonable to take in mind a slightly more generality:
we consider modules with highest weights over
 $\frg$,  whose restriction
 to $\frk$ is the  sum (not necessary direct)
of modules with highest weight.}
 Decompose $\frk$ into a sum of lowering, diagonal (Cartan)
 and raising (step up) subalgebras, $\frk=\frn_-+\frh+\frn_+$. Consider the
 ideal $U(\frg)\frn_+$ and the quotient $U(\frg)/ U(\frg)\frn_+$ over this ideal.
The reduction algebra (or
{\it Mickelsson algebra} $S(\frg,\frk)$) is defined as a subspace
 of $\frn_+$-highest  vectors in $U(\frg)/ U(\frg)\frn_+$.

It is easy to make sure that $S(\frg,\frk)$ is an algebra indeed and
that it acts on the space $V^{\frn_+}$ of $\frn_+$-highest  vectors of any
$\frg$-module $V$.
Next\footnote{assuming the conditions of the previous footnote are satisfied;},
$\frg$-module  $V$
can be reconstructed from $S(\frg,\frk)$-module $V^{\frn_+}$. \SS

The important  tool
 for the study of reduction algebras is the
 ''{\it {extremal projector}}'' of R.~M.~Asherova--Yu.~F.~Smirnov--
V.~N.~Tolstoy (1971). For a Lie algebra $\frs\frl_2$ with a standard basis
 $e$,
$h$, and $f$ it was obtained by  P.~O.~L\"owdin, \cite{Low} in 1964,
\begin{equation}
\label{kh1} P=\sum_{n=0}^\infty \frac{f^ne^n}{\prod_{j=1}^n (h+j+1)}
=\prod_{n=1}^\infty \left( 1-\frac{fe}{n(h+n+1)} \right)
\end{equation}
This expression defines an operator in any finite-dimensional
$\frs\frl_2$-module%
\footnote{It is defined in a more general setting, but this requires
an additional discussion},
   the operator  $P$ satisfies the relations
$$
e P=Pf=0, \qquad P^2=P
$$
In other words, $P$ is the projector to the subspace of
highest weight vectors.

 R.~M.~Asherova, Yu.~F.~Smirnov and V.~N.~Tolstoy found an analogous
 expression for arbitrary semisimple algebra  $\frk$. Namely
 for arbitrary positive root $\alpha$ we consider the related
$\frs\frl_2$ subalgebra with the  basis $e_\alpha$, $h_\alpha$, $f_\alpha$,
then write down the corresponding projector $P_\alpha$, and
take the product $\prod
P_\alpha$ over all positive roots in
a correct%
\footnote{so-called, normal,} order.

\SS


The formula for $P$ contains a division by elements of
the enveloping algebra $U(\frk)$,  i.e.,
 the element $P$ does
not belong to $U(\frk)$.

In this connection, Zhelobenko introduces
a certain extension of the algebra
 $U(\frk)$ (with partially allowed division) and identifies
 it with locally-finite endomorphisms
 of some universal module; after this
 an existence and uniqueness
 of the projection operator becomes a tautology.

\SS

The involvement of the projector force to localize%
\footnote{i.e., to allow a division partially}
 the algebra  $S(\frg,\frk)$ over expressions $h_\gamma+k$, where $h_\gamma$ is a coroot,
 and $k$ is an integer. In the localized algebra $Z(\frg,\frk)$,
there appear
 natural generators $z_g=Pg$, where $g$ runs the orthogonal complement
to  $\frk$ in $\frg$.
 Generators $z_g$ satisfy quadratic-linear relations with coefficients being
  rational expressions   in
 Cartan elements. The complete or a partial knowledge
 of these relations allows  to obtain various information about
 representations of the Lie algebra $\frg$.
 For instance, a description of the Mickelsson algebra for
the pair
  $\frg\frl_{n+1}\supset
\frg\frl_n$   gives eventually a possibility to obtain explicit formulas
for the action of generators of
$\frg\frl_n$ in Gel'fand-Tsetlin basis (\cite{AST2}), while the investigation
of a pair, related to
 a symmetric space is applied  in
the classification of certain
 series of Harish-Chandra modules  over
 real Lie groups.%
\footnote{For other applications see  \cite{Molev2}.}

\SS

One may "transpose" the defining properties of generators $z_g$ and get another
generators $z'_g$, which also satisfy quadratic-linear relations.
Zhelobenko constructs an operator  $Q$ such that $z'_g=Qg$. The operator $Q$
(known to specialists as the "Zhelobenko cocycle") also admits a factorization
over the roots and reminds
$P$ in $sl_2$ case:
\begin{equation}\label{kh2}Q_\alpha(x)=\sum_{n=0}^\infty \mathrm{Ad}(e_\alpha)^n
x \cdot f_\alpha^n \frac{(-1)^n}{\prod_{j=1}^n j(h_\alpha-j+1)}\, ,
\end{equation}
where 
$\mathrm{Ad}(e_\alpha) x:= [e_\alpha,x]$  is an operator of the adjoint action.
Later operators (\ref{kh2}) appeared in mathematical physics under the name of
dynamical Weyl group, see  V. Tarasov, A. Varchenkob \cite{TV}, G.~Felder
\cite{F}, P.~Etingof, A.~Varchenko \cite{EV}, see also  \cite{KN}.
 Even so, the general understanding of the very formula
(\ref{kh2})  is still unsatisfactory.

\SS

The discussed book of Zhelobenko contains quite a few inaccuracies; besides,
the last chapter is, to our point of view, unsuccessful. But as a whole, the book is an
impressive and very original work.



\SS

{\bf  5. The book "Principal structures of the representation theory". 2004.}
The work is very unusual by the style and structure. This is a text of
 "patchwork structure", which was not originally presumed for a
"systematic study" in the proper sense of this word.
The book consists of fragments; each of them, however, is quite readable and
 interesting.  From one side, they are expected (written) for the beginners,
 on the other hand, the professional mathematician could find there quite a number
  of interesting and unexpected.

   \SS

{\bf 6. Indecomposable representations of the Lorentz group (1958--1959).}
Infinite-dimensional nonunitary representations of the Lorentz group
$\SL(2,\C)$ (just as of any other semisimple group) are not completely reducible
 (i.e,, are not decomposed into a direct sum of irreducible representations).
 In particular, for exceptional values
of the parameters (see above) the representation
 of principal series of $\SL(2,\C)$ splits into a finite-dimensional representation and
an infinite-dimensional ``tail''%
\footnote{which  by itself is a point of the principal series}.
   Two notes \cite{Zh-Lorentz1}--\cite{Zh-Lorentz2}
    of Zhelobenko in "Doklady AN SSSR" are devoted to the problem of description
   of all indecomposable representations of the Lorentz group (with a finite
   composition series). Only  subrepresentations with the same infinitesimal
   character could be  linked.

Zhelobenko found a way to reduce this problem
to  a problem of finite-dimensional linear algebra.
Namely, one should classify collections of linear operators $d_+:V\to
W$, $d_-:W\to V$,  $\delta:W\to W$ such that the block matrices
$$
a:=\begin{pmatrix} 0&d_-\\ d_+&0 \end{pmatrix}, \qquad b:=\begin{pmatrix} 0&0\\
0&\delta \end{pmatrix}
$$
satisfy the relations:
$$
ab=ba, \qquad \text{$a$, $b$ are nilpotent}
$$
In its turn, this problem of linear algebra turned out to be quite curious, it was
solved by I.~M.~Gel'fand
 and V.~A.~Ponomarev \cite{GP}.
(1968).


An analogous  reduction to finite-dimensional linear algebra
problem was done for the groups $\SO(n,1)$ and $\SU(n,1)$, see \cite{Hor}.
 For arbitrary real semisimple groups,
 a general method (being far from effective algorithm)
 was suggested in \cite{BGG}. Probably these results are not ultimate.


\SS

{\bf 7. Gel'fand-Tsetlin bases.} Finite-dimensional irreducible
representations of semisimple groups admit a simple parametrization
(Cartan highest weight theorem),
however individual representations are complicated objects
and the description of their explicit implementations
 is a rather sophisticated problem.
An analogous problem appeared earlier for symmetric groups  $S_n$.
In 1931,  Alfred Young \cite{You} suggested a way
of construction of models of representations
with a help of successive  restriction
of a representation to subgroups $ S_n\supset S_{n-1}\supset S_{n-2}\supset\dots $.

This method was applied by  I.~M.~Gel'fand and M.~L.~Tsetlin  \cite{GTs} , 1950,
to classical groups  $\GL(n,\C)$ and $\OO(n,\C)$. Recall that restricting
an irreducible finite-dimensional representation of  $\GL(n,\C)$ to $\GL(n-1,\C)$,
we get a multiplicity free direct sum. Then we restrict
each of the obtained representations to  $\GL(n-2,\C)$  and so on.
We get finally a decomposition of the initial space into direct sum of one-dimensional
spaces. Choosing in each of them a vector (not canonically) we get a basis.
 This part of arguments is relatively trivial%
 \footnote{One can refer to  the Pieri formula.}.
The deep result is the explicit formulas for
the action of generators of the Lie algebra
 $\frg\frl(n,\C)$ in this basis. Gel'fand and Tsetlin announced
 these formulas in  two notes in "Doklady AN SSSR"

Later on different authors published several proofs, the earliest are the work
of Zhelobenko \cite{Zhe-spectral} and of
G.~Baird, L.~Biedenharn \cite{BaBi}. More important is that Zhelobenko
managed to develop here methods applicable for the solutions
of other spectral problems (and this was one of the starting points
for the books \cite{Zhe-compact}, \cite{Zhe-reductive}, see also the survey
 \cite{Zhe-GTs-2}).%
\footnote{Note also that Zhelobenko (1962, see  \cite{Zhe-spectral},
\cite{Zhe-compact}) obtained the spectrum of the reduction
of a finite-dimensional representation of
$\Sp(2n,\C)$ to $\Sp(2n-2,\C)$.  It is not multiplicity free,
but resembles the spectrum of the restriction of $\OO(k,\C)$ to $\OO(k-2,\C)$;
this enforced to think that Gel'fand-Tseitlin bases for $\Sp(2n,\C)$ exist.
\newline
\obmanka
 The question was many times discussed. see e.q. an approach
suggested by
  A.~A.~Kirillov and realized by
V.~V.~Stepin \cite{Sht}, 1986. The problem was finally solved by
 A.~I.~Molev
\cite{Molev}, 1999.%
\newline
\obmanka
 Notice also that the note of Zhelobenko on this subject in
Russian Math. Surveys 42 (1987) is mistaken.}

\SS


{\bf 8. Holomorphic families of intertwining operators and
 "operators of discrete symmetry".} Finally, we discuss the works of Zhelobenko
 on intertwining operators. Partially they are included
 into the book  "{\it Harmonic analysis on complex groups}",
but there is a later paper
 \cite{Zhe-discrete} on real groups.

\SS

Above (see n.3) we described a construction of
intertwining operators for representations of the principal series
of $\GL(n,\C)$; for any permutation $\sigma\in S_n$ of the parameters
 of a representation it gives a meromorphic operator-valued function
of the parameters of the representation%
\footnote{or, according to Bruhat \cite{Bruhat}, a function with values in
$G$-invariant  distributions on the product
of flag manifolds.}
with possible poles  at the points of reducibility.
 For instance, for the group $\GL(n,\C)$ this is a function
 on $p_1$, \dots, $p_n$.  The numbers $p_j-q_j$ are integers,
 i.e., we have
  $\Z^n$ of meromorphic functions.

The behavior of such functions at the points
of reducibility is an important and
complicated question, discussed in a plenty  of papers.
Zhelobenko introduces (  for complex groups)
 unobvious normalization, making
this family of operators holomorphic and nonvanishing everywhere.
In particular, this gives explicit formulas for intertwining operators in
reducibility points
as well.

\SS

However, these are not all intertwining operators for the principal series.
Sometimes (only in reducibility points)
 there are operators,
 which connect representations
with different collections of $\{p_j-q_j\}$%
\footnote{It seems that Naimark \cite{Nai}
was the first who  observed
 that in the study of
degenerations of principal series.}. Zhelobenko realized
 that
 an existence of such
 exceptional symmetries is
related to the natural action of the double of symmetric
 group $S_n\times S_n$ on the set of collections
$p_1$, \dots, $p_n$, $q_1$, \dots, $q_n$.
 He also presents explicit constructions of these operators.

 For  arbitrary real semisimple Lie algebra $\frg$ this is related to
 the action of the Weyl group of the complexification of $\frg_\C$ on the dual space
 to the Cartan subalgebra%
\footnote{ An heuristic  explanation for the reader slightly familiar with the subject:
eigenvalues of Laplace operators are constant on such orbits}.

\vspace{22pt}

As far as we remember Moscow of 1970-80, D.~P.~Zhelobenko looked
 rather  detached person.
But many people, mathematicians and physicists, students and experts,
 in the Soviet Union and
abroad, read the book "{\it Compact Lie groups}".
 To our mind, Zhelobenko's activity seemed to
influence upon the formation, interest and tastes
of a quite a number of the young at that
time mathematicians:
 A.~V.~Zelevinsky, G.~L.~Litvinov,
A.~I.~Molev, M.~L.~Nazarov, Yu.~A.~Neretin,
G.~I.~Olshanski, V.~N.~Tolstoy, S.~M.~Khoroshkin, I.~V.~Cherednik,%
\footnote{Of course, authors confirm this relatively to themselves.} though
none of them was Zhelobenko's student neither formally nor in fact.

\vspace{30pt}

Our paper was written  on advise of
 R.~S.~Ismagilov.\,
M.~S.~Al-Nator, M.~Duflo, Tamara Ivanovna Zhelobenko,  A.~V.~Zelevinsky,
G.~L.~Litvinov,  A.~S.~Kazarov,  M.~L.~Nazarov, V.~R.~Nigmatullin,
V.~F.~Molchanov, G.~I.~Olshanski, V.~N.~Tolstoy, N.~Ya.~Khelemski, V.~I.~Chilin
also helped us.

\end{document}